# Modelling hidden structure of signals in group data analysis with modified (Lr, 1) and block-term decompositions

Pavel Kharyuk[a,b,c,*], Ivan Oseledets[a,b]

[a] *Skolkovo Institute of Science and Technology, Skolkovo Innovation Center, Building 3, Moscow, Russia, 143026*
[b] *Institute of Numerical Mathematics of the Russian Academy of Sciences, Gubkina st. 8, Moscow, Russia, 119991*
[c] *Faculty of Computational Mathematics and Cybernetics, Lomonosov Moscow State University, Leninskiye Gory 1-52, GSP-1, Moscow, Russia, 119991*

## Abstract

This work is devoted to elaboration on the idea to use block term decomposition for group data analysis and to raise the possibility of modelling group activity with (Lr, 1) and Tucker blocks. A new generalization of block tensor decomposition was considered in application to group data analysis. Suggested approach was evaluated on multilabel classification task for a set of images. This contribution also reports results of investigation on clustering with proposed tensor models in comparison with known matrix models, namely common orthogonal basis extraction and group independent component analysis.

### Highlights

- Generalized block term decomposition with Tucker and (Lr, 1) blocks was proposed
- Models were modified to achieve purposes of group data analysis as feature extractors
- In classification, proposed approaches showed comparable with matrix methods quality

## 1. Introduction

Large number of studies are connected with analysis of datasets collected from various specific fields of science and technology. Under conditions of such studies it is often necessary to extract group information from data. Examples fitted into this group may appear in form of multisubject studies in cognitive, medical and other brain investigations (e.g., [1–3]) or in form of multisample examinations in chemical studies ([4–6]). From such data it is natural to expect it to have common (group) and individual information, and the main objective of group data analysis can be defined as reconstruction of these parts from given dataset. However, such formulation requires a word of clarification.

At first, the group *component* [7] analysis was considered on the premise that data structure permits us to reveal these components as vectors. Natural choice is to use multidimensional arrays as input data: properly tuned tensor decompositions capture components as inherent part of a method. In case of dealing with data as a set of matrices, matrix factorization techniques can be used. For example, in [8] extraction of common components proceeded via sequential matrix decompositions in three steps: individual dimensionality reduction with principal component analysis, constructing "essential" principal components (PCs) from concatenated individual PCs, and feature extraction with independent component analysis (ICA). This approach makes it possible to get common components shared by all objects or subjects. Another recently

*Corresponding author
*Email addresses:* kharyuk.pavel@gmail.com (Pavel Kharyuk), ivan.oseledets@gmail.com (Ivan Oseledets)

introduced methods, namely common orthogonal basis extraction (COBE) [9] and joint and individual variation explained (JIVE) [10], allow both individual and common parts to be extracted.

But in case of higher number of dimensions it is more realistic to decompose multidimensional data using tensor decompositions. In [3] authors used constrained canonical polyadic decomposition (CPD) to compute factor matrices for spatial, temporal and subject / session modes of multisubject/multisession fMRI data. However, separation of components into group and individual parts was not incorporated in that approach.

Having techniques for extraction of group information, one can use it as a part of some sort of processing framework. However, adequacy of extracted results should be understood in order to put such technique into any analysis pipeline. Here, a greater focus was placed on using proposed models in comparison with several of methods mentioned above for group data analysis in two machine learning tasks, namely classification and clustering [11]. Multiclass classification refers to algorithm that derives specific rule to label input data as representatives of one class from a set of classes. In clustering task, all input data need to be separated into several partitions named clusters. These tasks are coherent with group data analysis: group part contains information related to the class which can be utilized for classification task. And one can remove group information from the whole dataset to highlight differences between inputs expecting that it will increase performance of clustering algorithms.

Approach proposed in this work is based on linked-multiway component analysis concept [12; 13] carried out as constrained block-term decomposition [14–16]. In this approach components are to be extracted in pre-specified modes from input data structured as multidimensional array. In our previous works [17; 18] on group data analysis constrained (Lr, 1) block term decomposition was considered as a part of classification pipeline. The name "block term decomposition" implies that given tensor is represented as sum of several terms in a specific format. These blocks were imposed to be associated with individual and common group information. Whole decomposition can be found as solution of certain constrained optimization problem. In this work, more general framework with mixed-type blocks (both Tucker and (Lr, 1) blocks) as well as full (Lr, 1) model were introduced and evaluated in classification and clustering tasks.

## 2. Preliminary and basic concepts

### *2.1. Mathematical notation*

In this work the following notation was used. Multidimensional real-valued array $X \in \mathbb{R}^{n_1 \times \ldots \times n_d}$ is called **tensor**. An element with multi-index $(i_1, \ldots, i_d)$ of tensor $X$ is designated as $X[i_1, \ldots, i_d]$. Its $k$-th dimension is referred as $k$-th mode of $n_k$ size.

Arbitrary tensor can be reshaped into a column vector with vectorization operation. The correspondence between indexes of original and vectorized result is setting up by the $f(i_1, \ldots, i_d) = 1 + \sum_{k=1}^{d} \left[ (i_k - 1) \prod_{l=1}^{k-1} n_l \right]$ rule: $\operatorname{vec}\big(X\big)[f(i_1, \ldots, i_d)] = X[i_1, \ldots, i_d]$, $i_k = \overline{1, n_k}$, $k = \overline{1, d}$. If the result is a row vector, a circumflex mark (a hat) is used: $\widehat{\operatorname{vec}}(\cdot)$.

For $d$-dimensional tensor with $d > 1$ it is possible to reshape it into matrix form by joining axes. Such operation can be defined in various ways, here the leave-one-out approach is used where all axes except one (e.g., $k$-th) are to be joined. Let us call such operation a $k$-mode unfolding: $\operatorname{unfold}_k\big(X^{n_1 \times \ldots \times n_k \times \ldots n_d}\big) = \widehat{X}^{n_k \times \prod_{l \neq k} n_l}$. A contracted notation is also applicable, $X_{(k)} = \operatorname{unfold}_k(X)$.

Moreover, if there is an $l$-th mode such that its size $n_l$ can be factorized, $n_l = \prod_{p=1}^{P_l} n_{l,p}$, one may perform *tensorization* of input array. This operation designated as $\operatorname{Tens}(\cdot)$, and as it takes place, rule for indexes should be provided. Let us assume that default rule is a reshaping with natural mode range in column-major (Fortran) order.

Let $X \in \mathbb{R}^{n_1^{(x)} \times \ldots \times n_k \times \ldots n_{d_1}^{(x)}}$, $Y \in \mathbb{R}^{n_1^{(y)} \times \ldots \times n_k \times \ldots n_{d_2}^{(y)}}$ be $d_1$-dimensional and $d_2$-dimensional tensors, and let $n_k$ be an equal sized mode shared by both tensors $X$ and $Y$. If this holds, then *$k$-mode product* $X \times_k Y = Z$ between $X$, $Y$ is defined by the following element-wise formula:

$$Z[\dots, i_{k-1}, i_{k+1}, \dots, j_{k-1}, j_{k+1}, \dots] = \sum_{l=1}^{n_k} X[\dots, i_{k-1}, l, i_{k+1}, \dots] Y[\dots, j_{k-1}, l, j_{k+1}, \dots], \tag{1}$$

where $i_p = \overline{1, n_p^{(x)}}$, $j_q = \overline{1, n_q^{(y)}}$, $p = \overline{1, d_1}$, $q = \overline{1, d_2}$. In literature such operation is also known as contraction [19].

Frobenius norm for tensor $X$ may be written through the following expression: $\|X\|_F = \sqrt{\langle \mathrm{vec}(X), \mathrm{vec}(X) \rangle}$, where $\langle \cdot, \cdot \rangle$ is a scalar (inner) product between two vectors. The scalar product is also expandable for original tensors and is equivalent to performing consequent contractions across all modes of equal-sized input. If one takes contraction by all modes except those in $\Omega$ set, such operation may be designated in the following way: $\langle X, Y \rangle_{\{\alpha \in \Omega\}} = X \times_{k \in \{1, \dots, d\} \setminus \Omega} Y$. Braces placed in subscript of expression from here on mean that objects listed inside them are excluded from the operation.

In addition to inner product, outer product $X \circ Y = Z$ between two tensors $X$, $Y$ is defined as: $Z[i_1, \dots, i_{d_1}, j_1, \dots, j_{d_2}] = X[i_1, \dots, i_{d_1}] \cdot Y[j_1, \dots, j_{d_2}]$.

There are two more products it is necessary to specify. For two matrices $A$ and $B$ of sizes $(n_1, n_2)$ and $(m_1, m_2)$ Kronecker product is

$$A \otimes B = \begin{bmatrix} A[1,1]B & \dots & A[1,n_2]B \\ \vdots & \ddots & \vdots \\ A[n_1,1]B & \dots & A[n_1,n_2]B \end{bmatrix}. \tag{2}$$

If both $A$ and $B$ share an equal number of columns ($n_2 = m_2$ holds), Khatri-Rao product which is useful for CPD is defined:

$$A \odot B = \begin{bmatrix} A[:, i_1] \otimes B[:, j_1] & \dots & A[:, i_d] \otimes B[:, j_d] \end{bmatrix}. \tag{3}$$

For convenience all designations used in this work were listed in Table 1. For a detailed tutorial on basic operations and tensor decompositions the reader is referred to [20].

### *2.2. Component extraction*

It is a common assumption that each data instance from a given dataset has a certain structure. That structure possibly admits some kind of approximations that are suitable for using in practice to recover information structured in that way. One of the powerful concepts in that regard is source separation, where basis consisting of specific components is supposed to exist. Such a basis includes main information about data and is usually used as features for further analysis. In signal processing, source separation problem can be solved differently. One way is to use linear mixing model fitted by matrix factorization. Indeed, having measurements arranged in a matrix $X$, one can decompose it into matrix product with specified properties: $X = AB^T$, $X \in \mathbb{R}^{n_1 \times n_2}$, $A \in \mathbb{R}^{n_1 \times r}$, $B \in \mathbb{R}^{n_2 \times r}$. Different constraints like smoothness or sparseness (e.g., [21; 22]) may be incorporated here, but let us turn our attention to the following two types.

The first well-known decomposition arises from *principal component analysis* (PCA) [23] which frequently used for dimensionality reduction and construction of linear uncorrelated factors:

$$\begin{gathered} X \approx AB^T \\ s.t. \quad B^T B = I_r,\ B^T (X^T X) B = \mathrm{diag}(s_1, \dots, s_r) \end{gathered} \tag{4}$$

PCA is in great demand, however, uncorrelated variables are not necessary independent, and diagonalization of covariance matrix is sufficient only for Gaussian variables to share such a property. In other cases, higher order statistics are used to approximate statistical independence, on basis of which *independent component analysis* (ICA) [24] is conducted. By means of independence measure $f(\cdot)$ ICA problem may be formulated

in such a way:

$$\begin{aligned} &X \approx AB^T \\ s.t. \quad &f(A[:,1],\ldots,A[:,r]) \to \min \end{aligned} \tag{5}$$

In all these models, matrix $X$ is decomposed into a product between mixing and source matrices, and one may consider it as one-mode component extraction. Straightforward generalization of this concept is bimodal decomposition: $X \approx A_1 Y^T \approx A_1 G A_2^T = G \times_{n_1} A_1 \times_{n_2} A_2$, where columns of $A_1$ and $A_2$ keep sources for the first and the second modes respectively, and $G$ is a mixing core. Such representation is naturally expandable into higher dimensions for $X \in \mathbb{R}^{n_1 \times \cdots \times n_d}$ and known as *Tucker decomposition*:

$$\begin{aligned} X \approx \quad & |[G; A_1, \ldots, A_d]| = \\ = \quad & G \times_{n_1} A_1 \times_{n_2} \ldots \times_{n_d} A_d \\ G \in \mathbb{R}^{r_1 \times \cdots \times r_d},\ & A_i \in \mathbb{R}^{n_i \times r_i},\ i = \overline{1,d}\ . \end{aligned} \tag{6}$$

If core tensor is super-diagonal and has equal sizes of each dimension, one comes to *canonical polyadic decomposition* (CPD):

$$\begin{aligned} X \approx |[\Lambda; C_1, \ldots, C_d]| \\ \Lambda = \mathrm{Diag}(\lambda_1, \ldots, \lambda_r) \in \mathbb{R}^{r \times \cdots \times r},\ C_i \in \mathbb{R}^{n_i \times r},\ i = \overline{1,d}. \end{aligned} \tag{7}$$

Finally, each tensor may be decomposed into a sum of other tensors, where each term is approximated via its own decomposition: $X \approx T_1 + T_2 + \ldots + T_R$. If each term has either a Tucker or canonical polyadic structure, such decomposition is called *block term decomposition* (BTD) [14–16]. It should be noted that if each $T_i$ is a CP block, they are additionally constrained to have rank-1 factor-matrices for a subset of modes in each term. Block term decomposition as constrained CPD was originally named (Lr, Lr, 1)- decomposition due to the reason it was considered in 3-dimensional case. In this work, generalized designation for $N$-dimensional case, (Lr, 1) is kept.

Practice shows that more complicated tensor models make it possible to capture more relevant information from multidimensional data. For example, recently BTD was used for modelling epileptic seizures [25], and it was demonstrated that (Lr, 1)-decomposition outperforms usual CP decomposition in case where one of temporal, spectral or spatial stationarity assumptions was violated.

## 3. Block term decomposition and group data analysis

### *3.1. Data structure*

In this part the requirements and assumptions on structure of measurements are formulated. Firstly, let it be in a form of multidimensional array, for instance, multichannel measurements of electric activity of brain (electroencephalogram, EEG) acquired in different experimental conditions of trial-to-trial paradigm have essentially 4 dimensions - time, channel, condition, trial. The second requirement is consistency: each measurement should have equivalent number of counts for each dimension.

In many cases, one works with data organized as matrices (two-dimensional tensors). For example, resting-state EEG data collected from multichannel recording is naturally represented as matrix of "time-point" and "channel" axes. Monochrome images of objects collected from different angles for each of them can be stored as a set of matrices of "pixel" by "frame" shape. For any of similar datasets let us denote each data matrix as $X_i$. Such modes may be referred as "sample" and "channel".

One of the major objectives of group data analysis can be formulated as extracting common information from coupled dataset. The question is what may be called "common information". In this work, widely accepted point of view in signal processing is adhered to: the assumption that total activity is produced by transformations of several sources represented as vectors. Adaptation of such concept to group data analysis involves separating of sources into common, $S_{\text{com}}$, and individual, $S_i$.

**Assumption 1.** *Observations can be approximated by functional dependency from common and individual*

*sources*

$$X_i \approx F_i(S_{\text{com}}, S_i) \tag{8}$$

Subsequent step is to assume that common and individual sources are separable in a certain sense. Conventionally, it is thought that the separability is additive.
**Assumption 2.** *Observations are dependent on common and individual sources additively.*

$$X_i \approx f_i(S_{\text{com}}) + g_i(S_i) \tag{9}$$

Finally is should be clarified, which $f_i$, $g_i$ are suitable to use.
**Assumption 3.** *Linear mixing of sources is sufficient for approximation:*

$$f_i(S_{\text{com}}) = S_{\text{com}} B_{f_i}^T, \quad g_i(S_i) = S_i B_{g_i}^T \tag{10}$$

Putting it all together, one comes to the following model of data structure:

$$X_i \approx S_{\text{com}} B_{f_i}^T + S_i B_{g_i}^T. \tag{11}$$

This model brings up the question how to find such common and individual sources and their mixing matrices. In [8] the answer is supplied for common part by group ICA algorithm. Guoxu et al [9] proposed a way to extract both individual and common activity via simultaneous approximation of group part which is orthogonal to any individual one. In [10] group part is estimated via low-rank decomposition of $X = \begin{bmatrix} X_1 & \dots & X_N \end{bmatrix}$, and each individual part comes from subsequent low-rank decomposition of $X_i$ cleared of group part.

### *3.2. Proposed model*

Recall the matrix model with assumptions (8) and (9):

$$X_i \approx f_i(S_{\text{com}}) + g_i(S_i) \tag{12}$$

Suppose that each $X_i$, $i = \overline{1, N}$ has equal number of columns and rows. If that is true, one can concatenate matrices $X_i$ along the new, group axis. Let us denote the resulting 3-dimensional array as $X$. Adapted structure model is

$$X \approx \sum_{i=\overline{1,N}} \big(f_i(S_{\text{com}}) + g_i(S_i)\big) \circ e_i, \tag{13}$$

where $e_i$ is an i-th column of identity matrix $I_N$, $\circ$ is an outer product.

Further step to be taken is detailing of functions $f_i(A_{\text{com}})$ and $g_i(A_i)$. At first, let us focus on the individual part of sum, $\sum\limits_{i=\overline{1,N}} g_i(S_i) \circ e_i$. If one presumes linear mixing (assumption (10)), then $g_i(S_i) = S_i B_{g_i}^T$ holds, what could be directly captured in block-term model (namely, imposed (Lr, 1) decomposition). Imposed condition is a specific structure of 3-rd mode factor matrix, the identity one.

Considering group activity, the simplest assumption to be made is that $f_i(S_{\text{com}})$ differ only by (positive) constants $p_i$, $i = \overline{1, N}$:

$$f_i(S_{\text{com}}) = p_i f(S_{\text{com}}). \tag{14}$$

Putting linear mixing model here and collecting $p_i$ into column vector $p$, one comes to the following model:

$$X \approx S_{\text{com}} B^T \circ p + \sum_{i=\overline{1,N}} S_i B_{g_i}^T \circ e_i, \tag{15}$$

that can be rewritten as

$$\begin{aligned} X &\approx |[C_1, C_2, C_3]|, \\ C_1 &= \begin{bmatrix} S_1 & \dots & S_N & S_{\text{com}} \end{bmatrix}, \\ C_2 &= \begin{bmatrix} B_{g_1} & \dots & B_{g_N} & B \end{bmatrix}, \\ C_3 &= \begin{bmatrix} I_N & p \end{bmatrix} E, \end{aligned} \tag{16}$$

where $L = \begin{bmatrix} L_1 & \dots & L_N & L_{N+1} \end{bmatrix}$ specifies number of components for each individual and one common terms, and multiplication by $E = \begin{bmatrix} \mathbb{1}_{1\times L_1} \otimes e_1 & \dots & \mathbb{1}_{1\times L_{N+1}} \otimes e_{N+1} \end{bmatrix}$ replicates corresponding columns of $C_3$. This representation makes it clear how to generalize the model for higher-order inputs $X_i$:

$$\begin{aligned} X &\approx |[C_1, \dots, C_d]|, \\ C_i &= \begin{cases} \left[ \overbrace{C_i^{[1]}}^{L_1}, \dots, \overbrace{C_i^{[N+1]}}^{L_{N+1}} \right] & i = \overline{1, P} \\ \left[ c_j^{[1]}, \dots, c_j^{[N+1]} \right] E, & i = \overline{P+1, d-1} \\ & j = i - P \\ \begin{bmatrix} I_N & p \end{bmatrix} E, & i = d \end{cases} \end{aligned} \tag{17}$$

Here, the $d$-dimensional (Lr, 1) decomposition was used which was generalized in [26] with parameter $P$ corresponding to number of modes with full-sized factor matrices; the rest $d - P$ modes have reduced factor matrices.

Let us return to 3-dimensional case. If it is necessary to preserve mixing diversity, one should consider different mixing matrices for common sources leading to the $f_i(S_{\text{com}}) = S_{\text{com}} \hat{B}_i^T$. Stacking $\{\hat{B}_i\}_{i=\overline{1,N}}$ across group axis, one obtains a tensor $G$:

$$X \approx S_{\text{com}} \times_1 G + \sum_{i=\overline{1,N}} S_i B_{g_i}^T \circ e_i, \tag{18}$$

where $\times_1$ is a mode-1 multiplication. In such model there are $N$ (Lr, 1) blocks storing individual activity and one Tucker block for group information. In general, if block term decomposition contains $M$ Tucker terms and $N$ (Lr, 1) terms with common number of full-sized factor matrices, $P$, altogether, it is referred to as Tucker-(Lr, 1) decomposition (TLD) in this work:

$$\begin{aligned} X &\approx \sum_{m=\overline{1,M}} |[G^{[m]}; A_1^{[m]}, \dots, A_d^{[m]}|] + |[C_1, \dots, C_d]|, \\ C_i &= \begin{cases} \left[ \overbrace{C_i^{[1]}}^{L_1}, \dots, \overbrace{C_i^{[N]}}^{L_N} \right] & i = \overline{1, P} \\ \left[ c_j^{[1]}, \dots, c_j^{[N]} \right] E, & i = \overline{P+1, d} \\ & j = i - P \end{cases} \\ E &= \begin{bmatrix} \mathbb{1}_{1\times L_1} \otimes e_1 & \dots & \mathbb{1}_{1\times L_N} \otimes e_N \end{bmatrix}. \end{aligned} \tag{19}$$

As in full-(Lr, 1) block model, TLD-based model for group component analysis adopts the shape of the following representation:

$$X \approx |[G; A_1, \dots, A_d|] + |[C_1, \dots, C_d|],$$
$$C_i = \begin{cases} \Big[\overbrace{C_i^{[1]}}^{L_1}, \dots, \overbrace{C_i^{[N]}}^{L_N}\Big], & i = \overline{1, P}, \\ \Big[c_j^{[1]}, \dots, c_j^{[N]}\Big] E, & i = \overline{P+1, d}, \\ & j = i - P, \end{cases} \tag{20}$$
$$A_d = \mathrm{diag}(p),$$

where $p \in \mathbb{R}^{N \times 1}$ is a positive vector. In both models the last mode is reserved for group axis.

Note that although one can model individual parts with Tucker structure, this involves an increase in number of user-specified parameters (namely, Tucker ranks) for each term instead of vector $L$ for all individual terms.

In addition to hyperparameters considered above ($P$, $L$, and/or Tucker rank), it is necessary to select modes of interest across the first $P$ axes, $\Omega \subseteq \{1, \dots, P\}$. In order to separate common and individual information, additional constraints should be imposed on modes of interests. Such separation can be achieved in various ways. Following the lead of [9], the orthogonality constraint was explored:

$$S_{\mathrm{com};\gamma}^T S_{i;\gamma} = 0, \quad \forall i = \overline{1, N}, \forall \gamma \in \Omega. \tag{21}$$

### *3.3. Application in classification and clustering tasks*

In Figure 1, a sketch of general pipeline for using common and individual components in application to classification and clustering is presented. Classification task is seeking for the rule by which input data may be labelled as accurate as possible in terms of selected cost. Such rule is the result of training on a dataset with already known labels. If dataset is compatible with assumptions used in component models for group data analysis considered above, one can divide dataset by labels and extract group components for each part. After additional feature extraction in one mode of interest (that implies $\Omega = \{\gamma\}$) via ICA or other component extraction method, the first principle angle, which is related to canonical correlation analysis (CCA) [27], between column spaces of input and extracted features was used as distance measure:

$$\rho(Z, S_{\mathrm{com}}) = \min_{u,v} \left\{ \arccos \frac{|(u,v)|}{\|u\|\|v\|} \middle| u \in \mathcal{L}(\mathrm{columns}[Z]), v \in \mathcal{L}(\mathrm{columns}[S_{\mathrm{com}}]) \right\}. \tag{22}$$

$Z$ is a preprocessed instance from dataset, in this study it is a simple $n_\gamma$-mode unfolding, with mode of interest $\gamma$, $Z = \mathrm{unfold}_\gamma(Y)$.

Desired rule is preferring the class with the smallest principal angle ("winner takes all" strategy). In case of multiple modes of interest, principal angles for each of them may be combined using similar or some other strategy (for instance, voting).

In contrast to classification, clustering is an unsupervised learning problem with the objective to arrange the closest unlabeled elements into groups by means of certain affinity. It can be done differently, and one way is to use affinity matrices where one measures distances (or affinities) between each pair of objects to be clustered. Agglomerative clustering [11], which was used in this study, subsequently merges each two neighbouring clusters starting from distinct samples and up to joining all of them into one final cluster. At each step distances are recomputed according specified linkage rule.

Major problem in this task is appropriate selection of such linkage and affinity measure. In this study, several of them were considered: $l_1$, $l_2$, canberra, cosine, correlation and exponential map of $l_2$ distances; average and complete linkages.

### *3.4. Implementation*

All models and computational experiments were implemented using Python programming language under Anaconda distribution [28], which includes various pre-built packages for scientific computing. In this study the following packages were used: numpy [29], scipy [30], scikit-learn [31], matplotlib [32], seaborn [33]. All scripts can be found at GitHub repository: `https://github.com/kharyuk/gbtd`. Computational experiments are structured as Jupyter Notebooks [34].

## 4. Results and discussion

At first, implemented approaches were tested in experiments with random data. Tensor to be approximated was reconstructed from artificial random parameters. The following hyperparameters were used:

- Tensor configuration: $n_i = 20$, $i = \overline{1, d}$, $d = 3$;
- (Lr, 1) terms: $L_r = 3$, $r = \overline{1, 5}$, $P = d - 1 = 2$, (5 instances);
- Tucker term: $r_{T;i} = 3$ (1 instance).

Further experiments were done with group constraints. To conform to assumptions of group data analysis, the last mode of tensor was changed to $n_d = N = 5$ (5 objects/subjects) with consequent changes in $r_{T;d} = n_d$ for TLD and $L_{N+1} = n_d$ for full (Lr, 1) decomposition. Separation criterion was set on for the first mode, $n_1$.

Classification and clustering experiments were performed on ETH-80 database [35; 36], which includes color images separated into 8 categories with 10 objects per category depicted from different angles (41 frames per object). Each image was vectorized to form joint pixel mode. Therefore, information about each object was represented as $41 \times 128 \cdot 128 \times 3$ array (frame, pixel, color modes). Separation criterion was applied to the pixel mode.

### *4.1. Convergence on artificial data*

Optimization methods, implemented for models proposed in this study, are: alternating least squares (ALS), gradient descend (GD), conjugate gradients (CG) with 4 formulas for computing parameter $\beta$ (Fletcher-Reeves, CG-FR; Polak-Ribiere, CG-PR; Hestenes-Stiefel, CG-HS; Dai-Yuan, CG-DY), Gauss-Newton (GN), Levenberg-Marquardt method (quadratic, LM-Q; and Nielsen, LM-N), trust region with dogleg path (TR-DL), and Steihaug's conjugate gradients method with approximated (SCG-QN) and full (SCG-FN) Hessian. Expressions for computing Hessian and incorporating group constraints can be found in Appendix A and Appendix B respectively.

Convergence plots for unconstrained optimization are given in Figure 2 (medians of 50 runs). First-order methods had worse convergence for full (Lr, 1) model in comparison with ALS/second-order methods. Best results for both models were achieved by trust region method with dogleg trajectory. Another observation is that approximation of decomposition with mixed types of terms is computationally harder.

On Figures 3, 4 results for group constrained decompositions are shown. Strong oscillating behavior of convergence plots displayed on Figure 3 is a consequence of projection which perturbate the descending direction in general case. Similarly to unconstrained problem, ALS and second-order methods tend to have better estimation of parameters. Figure 4, (a) demonstrates that it makes good sense to use full Hessian to achieve better approximation, but it is relevant in cases where the solution should be as more accurate as possible. For both models optimized with Lagrange multipliers, SCG method achieved the best performance. As Figure 4, (b) shows, relatively close results may be obtained with dogleg method for TLD model.

However, when models were applied to the real data, it became critical to achieve reasonable level of time costs. In our experience with the ETH-80 data ALS was the fastest one. It was not observed that the more sophisticated methods provided better convergence on a few first iterations, and that was the reason ALS was selected for further evaluation of proposed models with ETH-80 dataset.

### 4.2. ETH-80: classification using group parts

In experiments with randomized data it was assumed that initial data tensor has a structure exactly as modelled. Such knowledge is not available for real data in most cases and so different combinations of hyperparameters should be additionally tested. In this study, multiple combinations of common and individual ranks were considered for proposed models, $r_{\text{com}} = r_{\text{c}} = \overline{1, 10}$, $r_{\text{ind}} = r_{\text{i}} = \overline{1, 10}$. Having components for pixel mode extracted with COBE, group (Lr, 1) (GLRO), group Tucker-(Lr, 1) decomposition (GTLD), independent components were additionally extracted from them via Fast ICA method [37]. 4-fold cross validation scheme was used to evaluate performance of approaches.

All methods showed comparable results, achieving higher than 90% score values (Table 2). It should be noted that modelling group activity as a Tucker term increased classification quality in comparison with full (Lr, 1) model. As it might be expected, higher values of common rank and lower individual ranks tended to have higher performance (Figure 5), except for the GICA models where individual rank affected number of individual principal components before merging them into common feature space for further extraction procedure.

### 4.3. ETH-80: clustering with individual parts

As in classification experiments, dataset has been splitted into 4 parts. Clustering was performed with use of each three parts followed by taking mean of the results. Three quality metrics were evaluated on partitioning with exact number of original classes: adjusted Rand index (accuracy equivalent corrected for chance), adjusted mutual information, Fowlkes-Mallows score (geometric mean of precision and recall). Clustering results are summarized in Table 3. As for GICA model, which was originally designed for estimating group activity only, contrasting was performed through subtraction backprojected group independent components. Matrix approaches as well as proposed tensor-based models achieved best performance with canberra distance and complete linkage. With such hyperparameters, all approaches substantially increase clustering quality (in comparison with results for $\text{Raw}_2$). However, best results for raw data, which were observed when correlation distance and average linkage were used ($\text{Raw}_1$), were getting closer to proposed tensor-based approaches. Matrix approaches showed the best performance, with leading COBE algorithm. It may be assumed that separation criteria applied to more than one mode can increase performance of GLRO/GTLD. Another possibility is to change separation criteria and / or increase number of iterations (only 10 iterations were applied).

## 5. Conclusions

New models for linked-multiway component analysis have been proposed as a constrained block term decomposition, as well as a generalization of block tensor decomposition with multiple Tucker and (Lr, 1) blocks. Computation of decompositions' parameters was formulated as a non-linear least squares task. Compared with group ICA and COBE-ICA models, block-term decomposition for group data analysis with constrained Tucker term proved itself to be appropriate for modelling group activity in classification task. Although clustering results were not as impressive, suggested approaches contain tunable parameters such as selection of various source modes and separation criteria, which potentially may be adjusted to increase the performance, and it will be considered as a venue for future work. Another drawback of proposed models is a higher computational complexity, which may be compensated for by the ability to extract mode specific factors.

## Acknowledgements

The work was supported by Russian Foundation for Basic Research (RFBR), research project No. 16-31-00494-mol_a. P.K. would also like to express gratitude to Dmitry Nazarenko for extensive proofreading.

## Appendix A: Structure of Hessian for Tucker-(Lr, 1) decomposition

Suppose that our task is to approximate given tensor $T \in \mathbb{R}^{n_1 \times \ldots \times n_d}$ in the Tucker-(Lr, 1) decomposition (19) format. Let us call $x$ a column vector composed of vectorized parameters of the model:

$$x = \begin{bmatrix} \ldots & \widehat{\text{vec}}(C_k^{[s]}) & \ldots & \widehat{\text{vec}}(A_k^{[m]}) & \ldots & \widehat{\text{vec}}(G^{[m]}) & \ldots \end{bmatrix}^T, \tag{23}$$

where $C_k^{[s]}$ is a $k$-th factor matrix within $s$-th CP term, $A_k^{[m]}$ is a $k$-th factor matrix within $m$-th Tucker term, $G^{[m]}$ is a Tucker core of $m$-th Tucker term, and let $F = F(x)$ be recovered tensor with current (approximated) values of parameters: $F(x) = \sum\limits_{m=\overline{1,M}} |[G^{[m]}; A_1^{[m]}, \ldots, A_d^{[m]}|] + |[C_1, \ldots, C_d|]$ (as in 19). For further conveniance let us use the following designation, $Z(x) = F(x) - T$. Corresponding optimization task may be formulated as non-linear least squares problem:

$$\min_x f(x), \quad f(x) = \frac{1}{2}\|F(x) - T\|_F^2 \tag{24}$$

Being the instance of nonlinear least squares task, unconstrained problem (24) has special expressions for gradient:

$$\nabla f(x) = J^T(x)\text{vec}(F(x) - T), \tag{25}$$

and Hessian matrix:

$$\begin{gathered} H[f](x) = J^T(x)J(x) + Q(x) \\ Q(x) = \sum_i \text{vec}(Z(x))[i] \cdot \nabla^2\left(\text{vec}(Z(x))[i]\right), \end{gathered} \tag{26}$$

where $J(x)$ is a Jacobi matrix for $f(x)$:

$$\begin{gathered} J(x) = \begin{bmatrix} \ldots & (V_{\{k\}}^{[s]} \otimes I_{n_k}) & \ldots & (G_{(k)}^{[m]} V_{\{k\}}^{[m]\,T} \otimes I_{n_k}) & \ldots & V^{[m]} & \ldots \end{bmatrix}, \\ k = \overline{1,d},\ m = \overline{1,M},\ s = \overline{1,S} \end{gathered} \tag{27}$$

with $V_{\{k\}}^{[s]} = C_d^{[s]} \odot \ldots \odot C_{k+1}^{[s]} \odot C_{k-1}^{[s]} \odot \ldots \odot C_1^{[s]}$, $V_{\{k\}}^{[m]} = A_d^{[m]} \otimes \ldots \otimes A_{k+1}^{[m]} \otimes A_{k-1}^{[m]} \otimes \ldots \otimes A_1^{[m]}$, $V^{[m]} = A_d^{[m]} \otimes \ldots \otimes A_1^{[m]}$.

The first part of Hessian has the following structure:

$$J^T J = \left[\begin{array}{c|c|c} \left(Gr^{CC}\right)_{k_1,k_2;s_1,s_2;\cdot,\cdot} & \left(Gr^{CA}\right)_{k_1,k_2;s_1,\cdot;\cdot,m_2} & \left(Gr^{CG}\right)_{k_1,\cdot;s_1,\cdot;\cdot,m_2} \\ \hline \left(Gr^{AC}\right)_{k_1,k_2;\cdot,s_2;m_1,\cdot} & \left(Gr^{AA}\right)_{k_1,k_2;\cdot,\cdot;m_1,m_2} & \left(Gr^{AG}\right)_{k_1,\cdot;\cdot,\cdot;m_1,m_2} \\ \hline \left(Gr^{GC}\right)_{k_1,k_2;\cdot,s_2;m_1,\cdot} & \left(Gr^{GA}\right)_{\cdot,k_2;\cdot,\cdot;m_1,m_2} & \left(Gr^{GG}\right)_{\cdot,\cdot;\cdot,\cdot;m_1,m_2} \end{array}\right], \tag{28}$$

$$k_1, k_2 = \overline{1, d},\ m_1, m_2 = \overline{1, M},\ s_1, s_2 = \overline{1, S}$$

$$Gr^{CC}_{k_1,k_2;s_1,s_2;\cdot,\cdot} = (V^{[s_1]\,T}_{\{k_1\}} \otimes I_{n_{k_1}}) P^{k_1,k_2}_{1,k_2} (V^{[s_2]}_{\{k_2\}} \otimes I_{n_{k_2}}), \tag{29}$$

$$Gr^{CA}_{k_1,k_2;s_1,\cdot;\cdot,m_2} = (V^{[s_1]\,T}_{\{k_1\}} \otimes I_{n_{k_1}}) P^{k_1,k_2}_{1,k_2} (V^{[m_2]}_{\{k_2\}} G^{[m_2]\,T}_{(k_2)} \otimes I_{n_{k_2}}), \tag{30}$$

$$Gr^{CG}_{k_1,k_2;s_1,\cdot;\cdot,m_2} = (V^{[s_1]\,T}_{\{k_1\}} \otimes I_{n_{k_1}}) P^{k_1}_{1} (V^{[m_2]}), \tag{31}$$

$$Gr^{AA}_{k_1,k_2;m_1,m_2} = (G^{[m_1]}_{(k_1)} V^{[m_1]\,T}_{\{k_1\}} \otimes I_{n_{k_1}}) P^{k_1,k_2}_{1,k_2} (V^{[m_2]}_{\{k_2\}} G^{[m_2]\,T}_{(k_2)} \otimes I_{n_{k_2}}), \tag{32}$$

$$Gr^{AG}_{k_1,\cdot;\cdot,\cdot;m_1,m_2} = (G^{[m_1]}_{(k_1)} V^{[m_1]\,T}_{\{k_1\}} \otimes I_{n_{k_1}}) P^{k_1}_{1} (V^{[m_2]}), \tag{33}$$

$$Gr^{GG}_{\cdot,\cdot;\cdot,\cdot;m_1,m_2} = (V^{[m_1]\,T})(V^{[m_2]}), \tag{34}$$

$P^{\cdot}_{\cdot}$ are commutation ("mode transition") matrices where upper and lower indices indicate modes and new places of these modes respectively. For example, vectorized transposition of matrix $A$ with modes may be written as $P^{1,2}_{2,1} vec(A) = vec(A^T)$.

In practical optimization $Q(x)$ term is usually omitted for several reasons: it takes additional computational expenses, selected computational method relies on symmetric positive definite approximation of Hessian, and because the closer solution attempts to local minimum, the less influence it has on the Hessian. Nevertheless, expressions to compute quadratic part of Hessian are also provided using tensor representation:

$$\text{Tens}\Big(\frac{\partial^2 z}{\partial \widehat{\text{vec}}(C^{[s_1]}_k) \partial \widehat{\text{vec}}(C^{[s_2]}_l)} \cdot z\Big) = \begin{cases} \langle |[C^{[s_1]}_1 \ldots, C^{[s_1]}_d]| |_{C^{[s_1]}_k, C^{[s_1]}_l = I_{L_{s_1}}}, Z \rangle_{\{k,l\}}, & k \neq l, s_1 = s_2 \\ 0, & \text{otherwise} \end{cases} \tag{35}$$

$$\text{Tens}\Big(\frac{\partial^2 z}{\partial \widehat{\text{vec}}(A^{[m_1]}_k) \partial \widehat{\text{vec}}(A^{[m_2]}_l)} \cdot z\Big) = \begin{cases} \langle |[G^{[m_1]}; A^{[m_1]}_1, \ldots, A^{[m_1]}_d]| \Big|_{\substack{A^{[m_1]}_k = I_{r^{[m_1]}_k}, \\ A^{[m_1]}_l = I_{r^{[m_1]}_l}}}, Z \rangle_{\{k,l\}}, & k \neq l,\ m_1 = m_2, \\ 0, & \text{otherwise} \end{cases} \tag{36}$$

$$\text{Tens}\Big(\frac{\partial^2 z}{\partial \widehat{\text{vec}}(A^{[m_1]}_k) \partial \widehat{\text{vec}}(G^{[m_2]})} \cdot z\Big) = \begin{cases} |[Z; A^{[m_1]}_1, \ldots, A^{[m_1]}_d]| \big|_{A^{[m_1]}_k = I_{r^{[m_1]}_k}}, & m_1 = m_2 \\ 0, & m \neq m_2 \end{cases} \tag{37}$$

where $m_1, m_2 = \overline{1, M}$, $k, l = \overline{1, d}$, $z = \text{vec}(Z(x))$, $Z = Z(x)$, and all other second derivatives are equal to zero.

It is worth to note that it may be very expensive to store Hessian matrix in computer memory. Indeed, to optimize memory costs it is reasonable not to work with the explicit representation but use a function which performs matrix-vector multiplication which suits any iterative solver, as used in our implementation.

## Appendix B: Incorporating group constraints into optimization task

In this appendix, it is described how to incorporate formulated constraints into optimization framework. Proposed model for group data analysis includes constraint on factor-matrices corresponding to group axis and separation constraint of the form (21). Let us define it explicitly for full-(Lr, 1) model (17):

$$\begin{gathered} C_d = \begin{bmatrix} I_N & p \end{bmatrix}, \\ p_1 + \ldots + p_N = p_{\text{cum}}, \quad p_i \geq p_{\min}, \\ \left(C_\gamma^{(N+1)}\right)^T C_\gamma^{(k)} = 0, \quad \gamma \in \Omega,\, k \in \overline{1, N}, \end{gathered} \tag{38}$$

and for TLD model (20):

$$\begin{gathered} A_d = \operatorname{diag}(p), \\ p_1 + \ldots + p_N = p_{\text{cum}}, \quad p_i \geq p_{\min}, \\ A_\gamma^T C_\gamma^{(k)} = 0, \quad \gamma \in \Omega,\, k \in \overline{1, N}, \end{gathered} \tag{39}$$

$p_{\min}$ and $p_{\text{cum}}$ are predefined constants, $\Omega$ is a set of modes of interest, $N$ is a total number of objects/subjects.

There are different ways to combine it with the problem (24), let us turn our attention to two methods: projected updates and Lagrange multipliers. In projected method, at each iteration current approximation is updated by solving the unconstrained problem followed by projection on feasible set defined by constraints. The first part of projector corresponding to separation constraint is a projection on the space orthogonal to linear span of subsequent matrix columns:

$$P_Y^{\perp}(X) = (I - Y(Y^T Y)^{-1} Y^T) X \tag{40}$$

The rest part is connected with projection onto the constrained $l_1$ ball [38] $\{y \in \mathbb{R}^N \,|\, y_1 + \ldots + y_N = p_{\text{cum}},\ y_i \geq p_{\min} > 0\}$.

Another way is to use Lagrange multipliers method where all constraints are incorporated as additional terms into the original unconstrained problem. Let us put all these terms into $g(x)$:

$$\begin{gathered} g(x) = \sum_{\gamma \in \Omega} \sum_{j=1}^{N} \left\langle vec(U_\gamma^T C_\gamma^{(j)}), \widehat{\mu}_{\gamma,j} \right\rangle + \left\langle \operatorname{vec}(\text{off diag}(Q_d)), \widehat{\lambda} \right\rangle - \\ -\theta \cdot \left(p_{\text{cum}} - \sum_{j=1}^{N} Q_d[j,j]\right) - \left\langle \zeta, \operatorname{diag}(Q_d) - p_{\min} \cdot \mathbb{1} \right\rangle, \\ U_\gamma = \begin{cases} C_\gamma^{(N+1)}, & \text{for model (17)} \\ A_\gamma, & \text{for model (20)} \end{cases}, \quad Q_d = \begin{cases} C_d, & \text{for model (17)} \\ A_d, & \text{for model (20)} \end{cases}, \end{gathered} \tag{41}$$

It is possible to reduce number of additional parameters, relaxing the constraints:

$$\begin{gathered} g(x) = \sum_{\gamma \in \Omega} \frac{\mu_\gamma}{2} \cdot \|U_\gamma^T C_\gamma\|_F^2 + \frac{\lambda}{2} \|\text{off diag}(Q_d)\|_F^2 - \\ -\theta \cdot \left(p_{\text{cum}} - \sum_{j=1}^{N} Q_d[j,j]\right) - \left\langle \zeta, \operatorname{diag}(Q_d) - p_{\min} \cdot \mathbb{1} \right\rangle, \end{gathered} \tag{42}$$

with the same $Q_d$, $U_\gamma$ as in (41). Gradient of constraint term for model (20) is defined by the following expressions:

$$\nabla g(x) = \begin{bmatrix} \nabla g_{\mu_1}(x) & \ldots & \nabla g_{\mu_{|\Omega|}}(x) & \nabla g_\lambda(x) & \nabla g_\theta(x) & \nabla g_{\zeta_1}(x) & \ldots & \nabla g_{\zeta_N}(x) \end{bmatrix} \tag{43}$$

$$\nabla g_{\mu_i}(x)^T = \begin{bmatrix} \ldots & 0 & \widehat{\operatorname{vec}}(A_i A_i^T C_i) & 0 & \ldots & 0 & \widehat{\operatorname{vec}}(C_i C_i^T A_i) & 0 & \ldots \end{bmatrix} \tag{44}$$

$$\nabla g_\lambda(x)^T = \begin{bmatrix} \ldots & 0 & \widehat{\operatorname{vec}}[\text{off diag}(A_d)] & 0 & \ldots \end{bmatrix} \tag{45}$$

$$\nabla g_\theta(x)^T = \begin{bmatrix} \ldots & 0 & \widehat{\operatorname{vec}}[\operatorname{diag}(\mathbb{1}_N)] & 0 & \ldots \end{bmatrix} \tag{46}$$

$$\nabla g_{\zeta_i}(x)^T = \begin{bmatrix} \dots & 0 & \underbrace{-1}_{i\text{-th position of } \mathrm{diag}(A_d)} & 0 & \dots \end{bmatrix} \tag{47}$$

Finally, the system with boarded Hessian is to be solved at each iteration:

$$\begin{bmatrix} H[f](x) & \nabla g(x) \\ \nabla g(x)^T & 0 \end{bmatrix} \begin{bmatrix} \Delta x \\ \Delta \tau \end{bmatrix} = \begin{bmatrix} -\nabla f(x) \\ -g(x) \end{bmatrix}, \tag{48}$$

where $\tau = \begin{bmatrix} \mu_1 & \dots & \mu_{|\Omega|} & \lambda & \theta & \zeta_1 & \dots & \zeta_N \end{bmatrix}^T$.

## Tables

Table 1: Notation used in the work

| Designation | Meaning |
|---|---|
| $\boxed{1}$ | array filled with ones |
| $(\dots)_{\{k\}}$ | operation with excluded $k$-th instance |
| $i = \overline{1, I}$ | integer variable $i$ takes values from $[1, I]$ |
| $A^+$ | Moore-Penrose pseudoinverse of matrix A |
| $\langle X, Y \rangle$ | inner product between two tensors (contraction) |
| $X \circ Y$ | outer product of two tensors |
| $A \otimes B$ | Kronecker product between two matrices |
| $A \odot B$ | Khatri-Rao product between two matrices |
| $X \times_k Y$ | $k$-mode tensor product |
| $X_{(k)}$ | $k$-mode unfolding of tensor $X$ |
| $\mathrm{vec}(X)$ | column vectorization of tensor $X$ |
| $\widehat{\mathrm{vec}}(X)$ | row vectorization of tensor $X$ |
| $\mathrm{diag}(v)$ | diagonal matrix with $v$ on main diagonal |
| $\mathrm{diag}(A)$ | main diagonal of matrix |
| off $\mathrm{diag}(A)$ | matrix A with zeros on main diagonal |
| $\|[C_1, \dots, C_d]\|$ | CP representation of $d$-dimensional tensor |
| $\|[G; A_1, \dots, A_d]\|$ | Tucker representation of $d$-dimensional tensor |

Table 2: Comparative results for best performing models in classification, mean of 4 fold CV. Precision, recall and F1 score were utilized for multiclass task with macro averaging.

| | COBE ($r_c = 9$) | GICA ($r_c = 9$, $r_i = 7$) | GLRO ($r_c = 9$, $r_i = 1$) | GTLD ($r_c = 10$, $r_i = 1$) |
|---|---|---|---|---|
| Accuracy | 0.943 | 0.948 | 0.911 | 0.938 |
| Precision | 0.958 | 0.964 | 0.930 | 0.945 |
| Recall | 0.943 | 0.948 | 0.911 | 0.938 |
| F1 score | 0.940 | 0.945 | 0.908 | 0.937 |

Table 3: Comparative results for best performing models on Agglomerative Clustering, mean of 4 runs on splitted data. COBE: $\mathrm{rank}_{\mathrm{com}} = 9$, canberra, complete; GICA: $\mathrm{rank}_{\mathrm{com}} = 10$, $\mathrm{rank}_{\mathrm{ind}} = 1$, canberra, complete; GTLD: $\mathrm{rank}_{\mathrm{com}} = 8$, $\mathrm{rank}_{\mathrm{ind}} = 3$, canberra, complete; GLRO: $\mathrm{rank}_{\mathrm{com}} = 8$, $\mathrm{rank}_{\mathrm{ind}} = 3$, canberra, complete; $\mathrm{Raw}_1$: correlation, average; $\mathrm{Raw}_1$: canberra, complete.

| | COBE | GICA | GTLD | GLRO | $\mathrm{Raw}_1$ | $\mathrm{Raw}_2$ |
|---|---|---|---|---|---|---|
| Adjusted Rand index | 0.693 | 0.615 | 0.583 | 0.565 | 0.545 | 0.423 |
| Adjusted Mutual Information | 0.775 | 0.728 | 0.693 | 0.669 | 0.669 | 0.556 |
| Fowlkes-Mallows score | 0.732 | 0.669 | 0.640 | 0.629 | 0.615 | 0.509 |

## Figures and captions

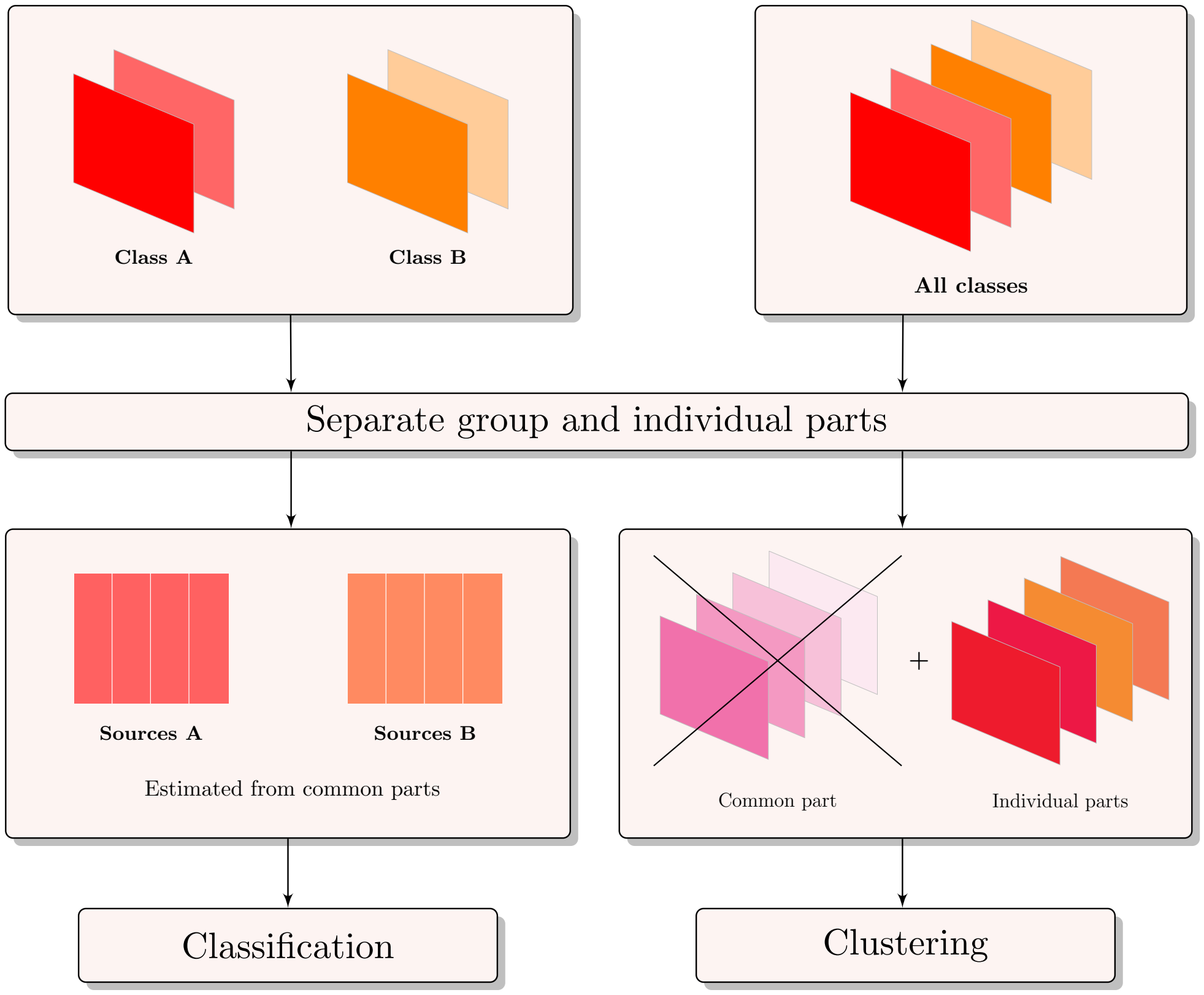

Figure 1: General classification and clustering pipeline based on source extraction and data contrasting.

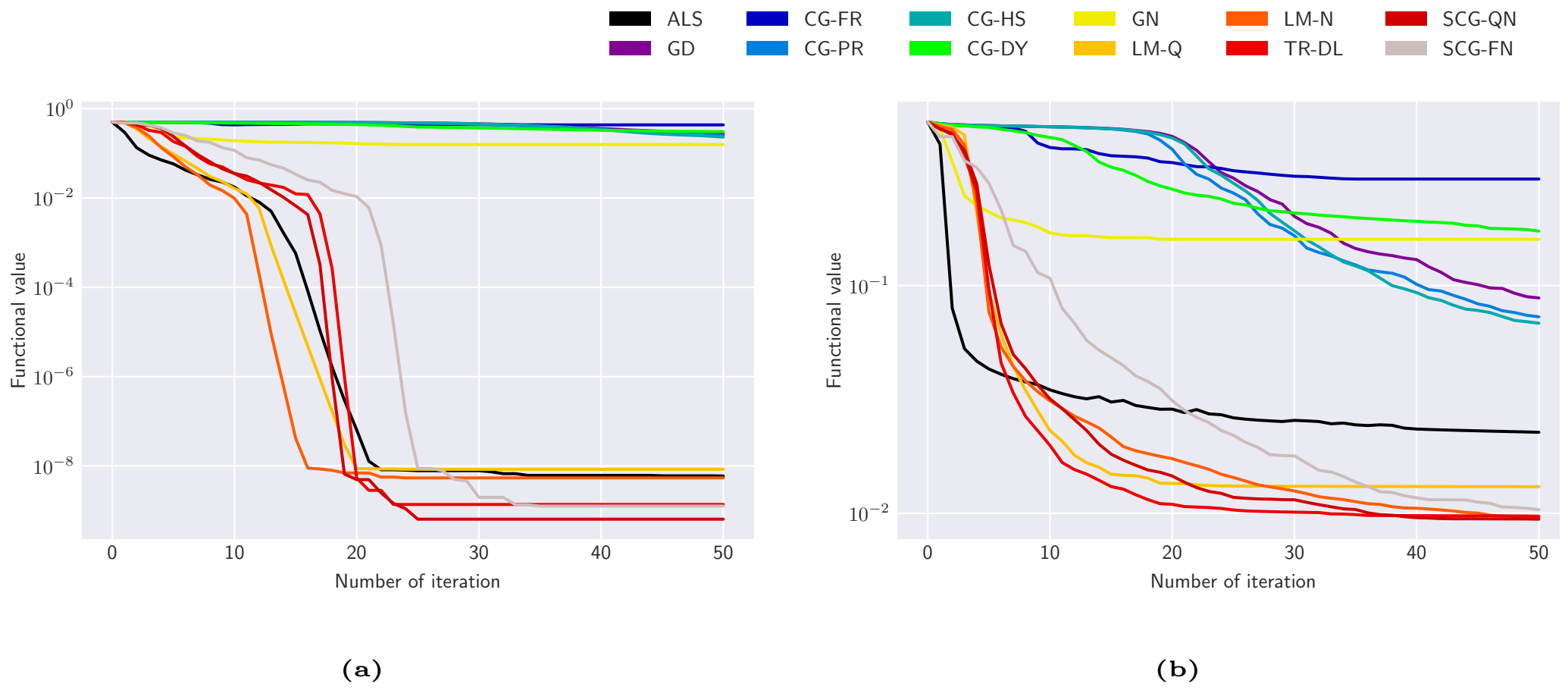

Figure 2: Convergence on artificial data, medians across 50 runs: (a) full (Lr, 1) decomposition; (b) Tucker-(Lr, 1) decomposition. Generating artificial data as well as parameters estimation were performed without constraints, using exact ranks.

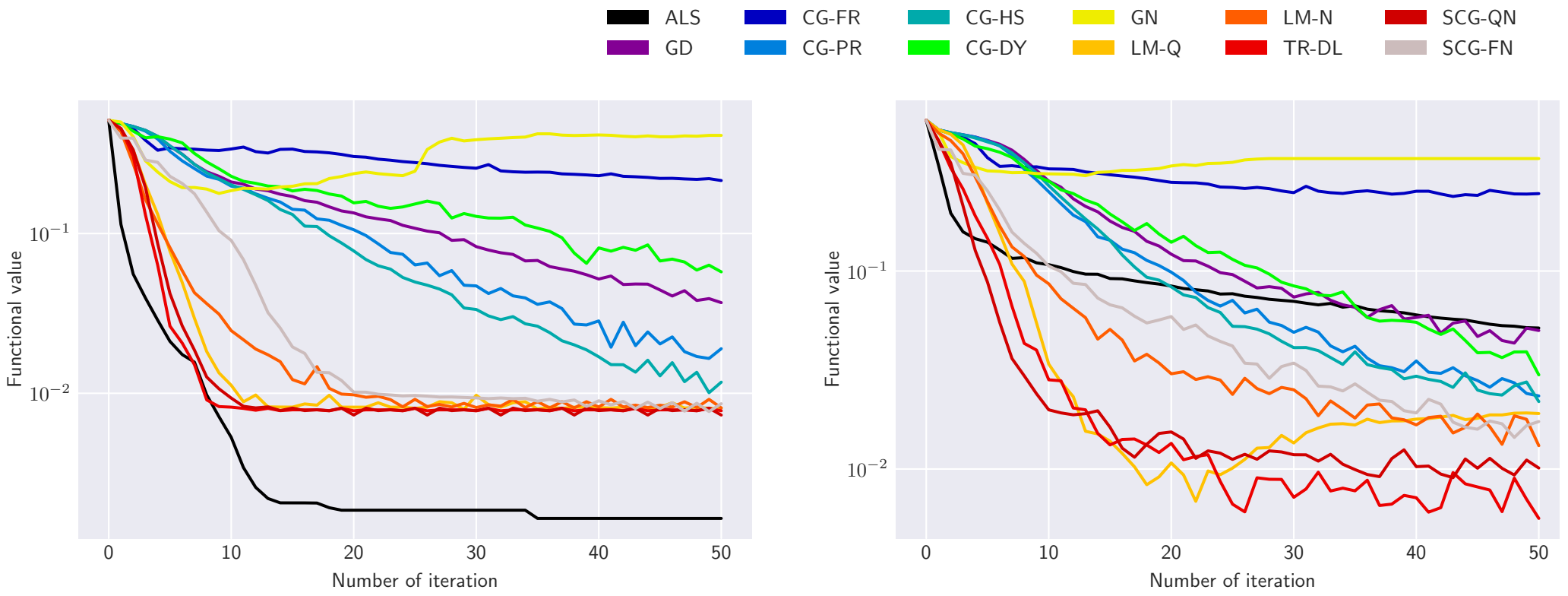

**(a)** **(b)**

Figure 3: Convergence on artificial data, medians across 50 runs: (a) full (Lr, 1) decomposition; (b) Tucker-(Lr, 1) decomposition. Generating artificial data as well as parameters estimation were performed with group constraints using exact ranks. Both decompositions were computed using projected iteration approach.

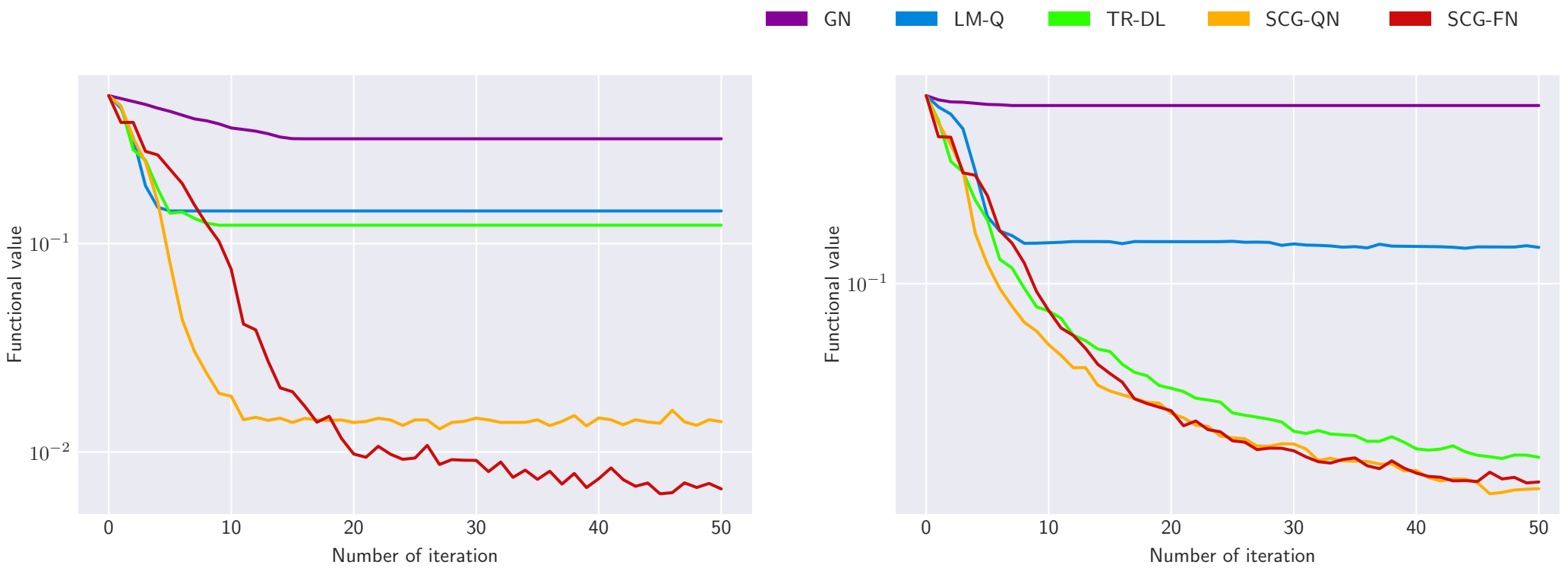

**(a)** **(b)**

Figure 4: Convergence on artificial data, medians across 50 runs: (a) full (Lr, 1) decomposition; (b) Tucker-(Lr, 1) decomposition. Generating artificial data as well as parameters estimation were performed with group constraints using exact ranks. Both decompositions were computed using Lagrange multipliers approach.

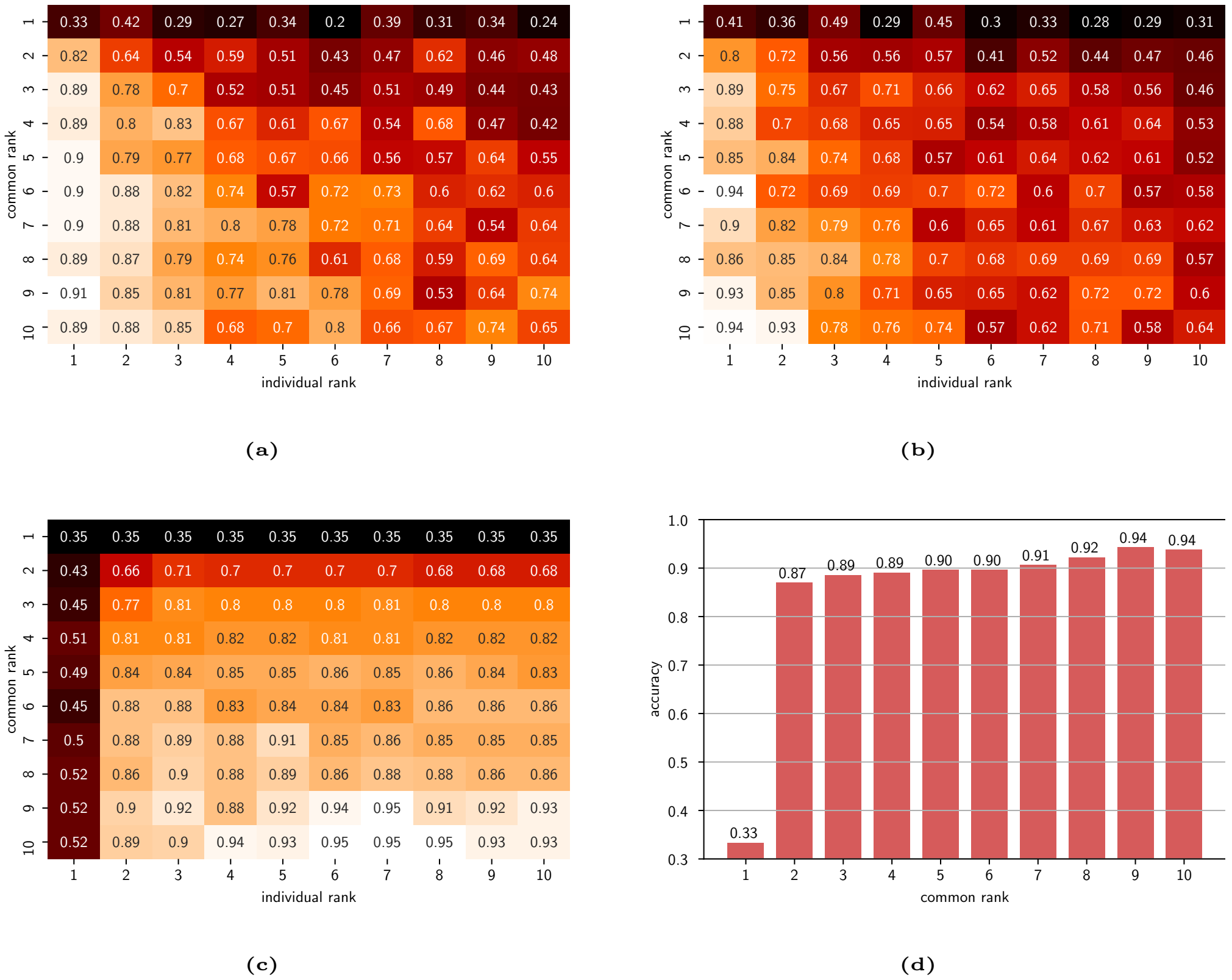

Figure 5: Mean performance (accuracy) of classification measured in conditions of 4-fold CV: (a) full (Lr, 1) decomposition; (b) Tucker-(Lr, 1) decomposition; (c) group ICA; (d) COBE.